\documentclass[a4paper]{amsart}
\usepackage{amsmath,amsthm,amssymb,latexsym,epic,bbm,comment,mathrsfs}
\usepackage{graphicx,enumerate,stmaryrd, color}
\usepackage[all,2cell]{xy}
\xyoption{2cell}

\newtheorem{theorem}{Theorem}

\newtheorem{proposition}[theorem]{Proposition}
\newtheorem{conjecture}[theorem]{Conjecture}

\usepackage[all]{xy}
\usepackage[active]{srcltx}
\usepackage[parfill]{parskip}
\usepackage{enumerate}

\begin{document}
\title[The tale of Kostant's problem]
{The tale of Kostant's problem}

\author[Volodymyr Mazorchuk]
{Volodymyr Mazorchuk}

\begin{abstract}
This is a survey paper presenting the history and both old 
and new results related to Kostant's problem. This problem
asks for which modules over a semi-simple finite dimensional 
complex Lie  algebra, the universal enveloping algebra surjects 
onto the algebra of adjointly locally finite linear endomorphism.
\end{abstract}

\maketitle

\section{Introduction, background and formulation}\label{s1}

\subsection{The general idea}\label{s1.1}

Let $A$ be a finitely generated associative algebra over the  complex numbers $\mathbb{C}$
and let $V$ be an $A$-module. The $A$-module structure on $V$ is given via the 
representation homomorphism $\varphi:A\to \mathrm{End}_\mathbb{C}(V)$.

If $V$ is finite dimensional and simple, the image of $\varphi$ is a simple
finite dimensional algebra over $\mathbb{C}$. Since the latter algebras are
exhausted by the matrix algebras due to Wedderburn--Artin theorem, it follows, in 
particular, that $\varphi$ is necessarily surjective.

If $V$ is infinite dimensional, then a simple dimension count shows that 
$\varphi$ is never surjective. Indeed, any finitely generated algebra has
either finite or countable dimension. At the same time, the dimension of the 
space of linear endomorphisms of an infinite dimensional complex vector 
space is always uncountable.

Consequently, it is natural to try to understand and describe the image of $\varphi$ 
inside $\mathrm{End}_\mathbb{C}(V)$ in some way. For example, one could try 
to single out this image using some axiomatic properties.

\subsection{The case of the universal enveloping algebra}\label{s1.2}

Let $\mathfrak{g}$ be a finite dimensional complex Lie algebra and
$U(\mathfrak{g})$ its universal enveloping algebra. Recall that the 
categories of $\mathfrak{g}$-modules and of $U(\mathfrak{g})$-modules are
naturally isomorphic. Furthermore, the associative algebra $U(\mathfrak{g})$
is finitely generated.

The algebra $U(\mathfrak{g})$ has the obvious structure of  a 
$U(\mathfrak{g})$-$U(\mathfrak{g})$-bimodule, given by multiplication.
In particular, $U(\mathfrak{g})$ is a $\mathfrak{g}$-$\mathfrak{g}$-bimodule.
Now, any $\mathfrak{g}$-$\mathfrak{g}$-bimodule $X$ can be considered as a
$\mathfrak{g}$-module with respect to the so-called {\em adjoint action}
$g(x):=g\cdot x-x\cdot g$, for $x\in X$ and $g\in\mathfrak{g}$.
We will denote this structure by $X^{\mathrm{ad}}$.
As a special case, the bimodule $U(\mathfrak{g})$ also has such a structure 
$U(\mathfrak{g})^{\mathrm{ad}}$ of an adjoint 
$\mathfrak{g}$-module.

From the commutation relations in $U(\mathfrak{g})$ it follows that 
the adjoint action of $\mathfrak{g}$ on $U(\mathfrak{g})$ is locally finite
(as it does not increase the degree of a monomial). This is a very striking
difference with both the left and the right regular actions 
of $\mathfrak{g}$ on $U(\mathfrak{g})$. 

Now let $M$ be a $\mathfrak{g}$-module. Then the vector space
\begin{displaymath}
\mathrm{End}_\mathbb{C}(M)=
\mathrm{Hom}_\mathbb{C}({}_{\mathfrak{g}}M_\mathbb{C},{}_{\mathfrak{g}}M_\mathbb{C})
\end{displaymath}
has the natural structure of a $\mathfrak{g}$-$\mathfrak{g}$-bimodule.
We also have that the 
representation map $\varphi:U(\mathfrak{g})\to \mathrm{End}_\mathbb{C}(M)$
is a homomorphism of $\mathfrak{g}$-$\mathfrak{g}$-bimodules.

The $\mathfrak{g}$-$\mathfrak{g}$-bimodule $\mathrm{End}_\mathbb{C}(M)$
has a subbimodule, usually denoted by $\mathcal{L}(M,M)$, which consists of
all elements of $\mathrm{End}_\mathbb{C}(M)$, the adjoint action of 
$\mathfrak{g}$ on which is locally finite (it follows from the Leibniz
rule that this linear subspace of $\mathrm{End}_\mathbb{C}(M)$ is, indeed, a
subbimodule, it is even a subalgebra). The image of $\varphi$ is, of course, a 
$\mathfrak{g}$-$\mathfrak{g}$-subbimodule of $\mathcal{L}(M,M)$.

\subsection{Kostant's problem}\label{s1.3}

In the above setup, Kostant's problem for a $\mathfrak{g}$-module $M$,
as formulated and popularized by Joseph in \cite{Jo}, is stated as follows:

{\bf Kostant's Problem.} Is the map 
$\varphi:U(\mathfrak{g})\to \mathcal{L}(M,M)$ surjective?

In other words, it asks whether the axiomatic description via adjointly finite
part is enough to single out the image of $\varphi$ from the rest of 
$\mathrm{End}_\mathbb{C}(M)$. Factoring the annihilator
$\mathrm{Ann}_{U(\mathfrak{g})}(M)$ of $M$ in $U(\mathfrak{g})$ out, 
one can equivalently ask whether the map $\varphi$ induces an isomorphism 
between  $U(\mathfrak{g})/\mathrm{Ann}_{U(\mathfrak{g})}(M)$
and $\mathcal{L}(M,M)$.

\subsection{Tools: Harish-Chandra bimodules}\label{s1.4}

In order to be able to answer the question  posed by Kostant's problem for 
some $\mathfrak{g}$-module $M$, we need some way to ``measure'' 
and compare the sizes of $U(\mathfrak{g})/\mathrm{Ann}_{U(\mathfrak{g})}(M)$
and $\mathcal{L}(M,M)$. This brings us to the notion of a 
Harish-Chandra $\mathfrak{g}$-$\mathfrak{g}$-bimodule.

A $\mathfrak{g}$-$\mathfrak{g}$-bimodule $X$ is called a
{\em Harish-Chandra bimodule} provided that 
\begin{itemize}
\item it is finitely generated;
\item the adjoint action of $\mathfrak{g}$ on $X$ is locally finite;
\item each simple finite dimensional $\mathfrak{g}$-module appears
in $X^{\mathrm{ad}}$ with finite multiplicity.
\end{itemize}

There is, of course, no guarantee that either of 
$U(\mathfrak{g})/\mathrm{Ann}_{U(\mathfrak{g})}(M)$
or $\mathcal{L}(M,M)$ is a Harish-Chandra bimodule, however, if
$\mathcal{L}(M,M)$ is a Harish-Chandra bimodule, then so is
$U(\mathfrak{g})/\mathrm{Ann}_{U(\mathfrak{g})}(M)$ as well,
and to answer Kostant's problem for such $M$ it is enough to 
compare the multiplicities of each simple finite dimensional 
$\mathfrak{g}$-module in the two adjoint $\mathfrak{g}$-modules
$\big(U(\mathfrak{g})/\mathrm{Ann}_{U(\mathfrak{g})}(M)\big)^{\mathrm{ad}}$
and $\mathcal{L}(M,M)^{\mathrm{ad}}$.

We note that $U(\mathfrak{g})$ itself is not a Harish-Chandra bimodule
since, for example, it has infinite dimensional center and hence the
multiplicity of the trivial $\mathfrak{g}$-module in the adjoint 
$\mathfrak{g}$-module $U(\mathfrak{g})$ is infinite. Therefore 
the use of Harish-Chandra bimodules in the study of Kostant's problem
is only reasonable for modules with fairly large annihilator, namely,
we need $\mathrm{Ann}_{U(\mathfrak{g})}(M)$ to be ``large enough''
to ensure that all simple finite dimensional $\mathfrak{g}$-modules
appear in 
$\big(U(\mathfrak{g})/\mathrm{Ann}_{U(\mathfrak{g})}(M)\big)^{\mathrm{ad}}$
with finite multiplicities.

\subsection{Kostant's Theorem}\label{s1.5}

A more precise content for what it could mean for the annihilator of 
$M$ to be ``large enough'' can be deduced from the following theorem
of Kostant:

\begin{theorem}[Kostant, \cite{Ko}]\label{kostantthm}
Assume that $\mathfrak{g}$ is semi-simple. 
Then the algebra $U(\mathfrak{g})$ is free as the module over the center
$Z(\mathfrak{g})$ of $U(\mathfrak{g})$. 
Moreover, for any simple finite dimensional $\mathfrak{g}$-module
$V$, the $Z(\mathfrak{g})$-module 
$\mathrm{Hom}_{\mathfrak{g}}(V,U(\mathfrak{g})^{\mathrm{ad}})$
is free and its rank 
equals the dimension of
the $0$-weight space $V_0$ in $V$ (in particular, it is finite).
\end{theorem}

Consequently, if $\mathrm{Ann}_{U(\mathfrak{g})}(M)\cap Z(\mathfrak{g})$
has finite codimension in $Z(\mathfrak{g})$, then the 
bimodule
$U(\mathfrak{g})/\mathrm{Ann}_{U(\mathfrak{g})}(M)$
is a Harish-Chandra bimodule. In particular, this means that
the approach to Kostant's problem via Harish-Chandra bimodules is 
potentially
applicable for all simple $\mathfrak{g}$-modules $M$ due to the fact
that all simple $\mathfrak{g}$-modules admit a central character,
see \cite[Proposition~2.6.8]{Di}.

\subsection{Canonical isomorphisms}\label{s1.6}

Assuming that $\mathfrak{g}$ is semi-simple and $V$ is a simple
finite dimensional $\mathfrak{g}$-module,
to measure the size of $\mathcal{L}(M,M)$, one can use the following
canonical isomorphisms, see \cite[Subsection~6.8]{Ja}:
\begin{equation}\label{eq2}
\mathrm{Hom}_\mathfrak{g}(V,\mathcal{L}(M,M)^\mathrm{ad})\cong
\mathrm{Hom}_\mathfrak{g}(M\otimes_\mathbb{C} V,M)\cong
\mathrm{Hom}_\mathfrak{g}(M,M\otimes_\mathbb{C} V^*).
\end{equation}
Here the tensor product of $\mathfrak{g}$-modules is considered as a 
$\mathfrak{g}$-module in the usual way and $V^*$ denotes the module
dual to $V$.

\section{Category $\mathcal{O}$ and Verma modules}\label{s2}

\subsection{Category $\mathcal{O}$}\label{s2.1}

From now on we assume that $\mathfrak{g}$ is semi-simple with a fixed
triangular decomposition
\begin{displaymath}
\mathfrak{g}=\mathfrak{n}_-\oplus \mathfrak{h}\oplus \mathfrak{n}_+, 
\end{displaymath}
where $\mathfrak{h}$ is a Cartan subalgebra of $\mathfrak{g}$.
Associated to this triangular decomposition we have the 
Bernstein-Gelfand-Gelfand (BGG) category $\mathcal{O}$ for
$\mathfrak{g}$, see \cite{BGG}, which is defined as the full
subcategory of the category of all $\mathfrak{g}$-modules that
consists of all $\mathfrak{g}$-modules $M$ such that
\begin{itemize}
\item $M$ is finitely generated;
\item the action of $\mathfrak{h}$ on $M$ is diagonalizable;
\item the action of $U(\mathfrak{n}_+)$ on $M$ is locally finite.
\end{itemize}

For a $\mathfrak{g}$-module $M$ and $\lambda\in\mathfrak{h}^*$ set
\begin{displaymath}
M_\lambda:=\{m\in M\,:\, h\cdot m=\lambda(h)m,\text{ for all }h\in\mathfrak{h}\}. 
\end{displaymath}
Such an element $\lambda$ is usually called a {\em weight} and the corresponding
$M_\lambda$ a {\em weight space}. If $M$ is the direct sum of its weight spaces,
then such $M$ is called a {\em weight module}. In other words, a weight module
is exactly a module the action of $\mathfrak{h}$ on which  is diagonalizable.
For a weight module $M$, the set of all weights of $M$ for which the corresponding
weight spaces are non-zero is called the {\em support} of $M$ and denoted
$\mathrm{supp}(M)$.

Let $\mathbf{R}$ be the root system of the pair $(\mathfrak{g},\mathfrak{h})$.
The triangular decomposition above induces a partition of 
$\mathbf{R}$ into the positive and negative roots:
$\mathbf{R}=\mathbf{R}_+\coprod \mathbf{R}_-$.

For a weight module $M$, a weight $\lambda\in \mathrm{supp}(M)$ is called
a {\em highest weight} provided that $\lambda+\alpha\not\in \mathrm{supp}(M)$,
for any $\alpha\in \mathbf{R}_+$. Any element in $M_\lambda$, for a highest weight
$\lambda$, is called a {\em highest weight vector}. The module $M$ is called a 
{\em highest weight} module provided that $M$ is generated by some highest weight vector.

\subsection{Verma modules}\label{s2.2}

Prominent objects of category $\mathcal{O}$ are Verma modules, see \cite{Ve}
and \cite[Section~7]{Di}.
For $\lambda\in\mathfrak{h}^*$, denote by $\mathbb{C}_\lambda$ the
one-dimensional $\mathfrak{h}$-module on which $\mathfrak{h}$ acts via $\lambda$.
This extends to an $\mathfrak{h}\oplus \mathfrak{n}_+$-module via
$\mathfrak{n}_+\cdot \mathbb{C}_\lambda=0$. The 
corresponding induced module
\begin{displaymath}
\Delta(\lambda):=
U(\mathfrak{g})\bigotimes_{U(\mathfrak{h}\oplus \mathfrak{n}_+)}\mathbb{C}_\lambda
\end{displaymath}
is called the {\em Verma module} with {\em highest weight} $\lambda$.
Indeed, $\Delta(\lambda)$ is a highest weight module and $\lambda$ is its unique
highest weight. Moreover, by adjunction, any highest weight module with 
highest weight $\lambda$ is a quotient of $\Delta(\lambda)$
(the so-called {\em universal property} of Verma modules).

The module $\Delta(\lambda)$ has a unique maximal submodule and hence a unique
simple quotient, usually denoted $L(\lambda)$. The set
$\{L(\lambda)\,:\,\lambda\in\mathfrak{h}^*\}$ is a complete and irredundant
set of representatives of isomorphism classes of simple modules in $\mathcal{O}$.

\subsection{Joseph's Theorem}\label{s2.3}

The first major result on Kostant's problem is the following theorem
of Joseph, see \cite[Corollary~6.4]{Jo}, which asserts that the solution to 
Kostant's problem is positive for all Verma modules:

\begin{theorem}[Joseph, \cite{Jo}]\label{josephthm}
For any $\lambda\in \mathfrak{h}^*$, the representation map induces an
isomorphism
\begin{displaymath}
U(\mathfrak{g})/\mathrm{Ann}_{U(\mathfrak{g})}(\Delta(\lambda))
\cong
\mathcal{L}(\Delta(\lambda),\Delta(\lambda)).
\end{displaymath}
\end{theorem}

To prove this theorem, one has to show that, for any simple finite dimensional
$\mathfrak{g}$-module $V$, the multiplicity of $V$ in 
$\big(U(\mathfrak{g})/\mathrm{Ann}_{U(\mathfrak{g})}(\Delta(\lambda))\big)^{\mathrm{ad}}$ 
and $\mathcal{L}(\Delta(\lambda),\Delta(\lambda))^{\mathrm{ad}}$ coincide.
In fact, as it turns out, both are equal to the dimension of $V_0$.

For the former bimodule 
$U(\mathfrak{g})/\mathrm{Ann}_{U(\mathfrak{g})}(\Delta(\lambda))$, this follows
by combining Kostant's Theorem (i.e. Theorem~\ref{kostantthm}) with a result of Duflo
(see \cite{Du} and \cite[Theorem~8.4.3]{Di}) 
that asserts that the annihilator of any Verma module is 
generated by the intersection of this annihilator with the center
of the universal enveloping algebra.

For the latter bimodule $\mathcal{L}(\Delta(\lambda),\Delta(\lambda))$,
the canonical isomorphisms from Subsection~\ref{s1.5} reformulate the
problem as
\begin{equation}\label{eq1}
\dim \mathrm{Hom}_\mathfrak{g}(\Delta(\lambda),\Delta(\lambda)\otimes_\mathbb{C} V^*) =
\dim V_0.
\end{equation}
Note that $\dim V_0=\dim V^*_0$. 
Understanding this formula requires a slightly deeper dive into 
the properties of category $\mathcal{O}$, so we will come
back to this a bit later.

\subsection{Block decomposition}\label{s2.4}

As already mentioned above, category $\mathcal{O}$ is quite big.
For instance, the isomorphism classes of simple objects in $\mathcal{O}$
are in bijection with the elements in $\mathfrak{h}^*$. This can be 
improved by decomposing $\mathcal{O}$ into direct summands. For this
we need to look more carefully into the action of $Z(\mathfrak{g})$.

Let $W$ be the Weyl group of $\mathbf{R}$. It acts naturally on
$\mathfrak{h}^*$, by definition. Let $\rho$ be the half of the sum
of all positive roots. The {\em dot-action} of $W$ on $\mathfrak{h}^*$
is defined as $w\cdot \lambda:=w(\lambda+\rho)-\rho$. A theorem by
Harish-Chandra, see \cite{HC} and 
\cite[Propositions~7.4.7 and 7.4.8]{Di}, asserts that
\begin{itemize}
\item any central character is realized as the central character of
some Verma module;
\item two Verma modules have the same central character if and only
if their highest weights belong to the same orbit of the dot-action of
$W$ on $\mathfrak{h}^*$.
\end{itemize}

For any central character $\chi:Z(\mathfrak{g})\to\mathbb{C}$, denote by
$\mathcal{O}_\chi$ the Serre subcategory of $\mathcal{O}$ generated by
all simple highest weight modules that admit this central character
$\lambda$. Then we have the decomposition
\begin{displaymath}
\mathcal{O}\cong\bigoplus_{\chi} \mathcal{O}_\chi,
\end{displaymath}
moreover, each $\mathcal{O}_\chi$ is equivalent to the category of 
finite dimensional modules over some finite dimensional associative 
$\mathbb{C}$-algebra, see \cite{BGG}. Note that some $\mathcal{O}_\chi$
might decompose further.

The set $\mathfrak{h}^*$ has a natural partial order: $\lambda\leq \mu$
if and only if $\mu-\lambda$ can be written as a linear combination of
positive roots with non-negative integer coefficients. Fix now some
$\mu\in \mathfrak{h}^*$ and consider the dot-orbit $W\cdot \mu$.
If $\lambda\in W\cdot \mu$ is a maximal element in this orbit with 
respect to the partial order above, then $\Delta(\lambda)$ is projective
in $\mathcal{O}$. 

As was already observed in \cite{BGG}, for any finite dimensional
$\mathfrak{g}$-module $V$, the module 
$\Delta(\lambda)\otimes_\mathbb{C}V$ has a finite filtration whose
subquotients are isomorphic to Verma modules. The corresponding 
multiplicities of Verma modules in $\Delta(\lambda)\otimes_\mathbb{C}V$
are well-defined (i.e. do not depend on the choice of such filtration).
Furthermore, the multiplicity of $\Delta(\lambda)$ equals
$\dim V_0$. This implies Formula~\eqref{eq1} (and hence Theorem~\ref{josephthm}
as well) in the case $\lambda\in W\cdot \mu$ is a maximal element.
The general case will need to wait a bit more.

\subsection{Quotients of projective Verma modules}\label{s2.5}

Another interesting application of projective Verma modules is the following
observation which can be found in \cite[Subsection~6.9]{Ja}, with a detailed
proof by direct computation: if
$\lambda\in \mathfrak{h}^*$ is such that $\Delta(\lambda)$ is projective
in $\mathcal{O}$, then, for any submodule $N\subset \Delta(\lambda)$,
the natural inclusion
\begin{displaymath}
U(\mathfrak{g})/\mathrm{Ann}_{U(\mathfrak{g})}(\Delta(\lambda)/N)
\subset
\mathcal{L}(\Delta(\lambda)/N,\Delta(\lambda)/N)
\end{displaymath}
is an isomorphism.

From this, it is very tempting to expect that the positivity of the 
answer to Kostant's problem for some module $M$ should be inherited
by all quotients of $M$. Unfortunately, this is wrong, in general.

\subsection{Principal block}\label{s2.6}

We have already decomposed $\mathcal{O}$ into direct 
summands $\mathcal{O}_\chi$ which are ``small enough'' to be
module categories of some finite dimensional algebras. These
can be further decomposed into indecomposable direct summands,
usually called {\em blocks}. For a fixed $\mathfrak{g}$, 
the corresponding category $\mathcal{O}$ has only finitely
many indecomposable direct summands, up to equivalence.

The ``most complicated'' one among these blocks is the so-called 
{\em principal block} $\mathcal{O}_0$ which is defined as the
indecomposable direct summand of $\mathcal{O}$ that contains
the trivial (one-dimensional) $\mathfrak{g}$-module. This 
trivial $\mathfrak{g}$-module is, in fact, the module $L(0)$,
where $0\in\mathfrak{h}^*$ is the trivial functional on $\mathfrak{h}$.

Isomorphism classes of simple modules in $\mathcal{O}_0$ are in bijection
with the elements of $W$, namely, they are $L(w\cdot 0)$,
for $w\in W$. As we will mainly concentrate on the principal 
block from now on, for simplicity, we will write 
$\Delta_w$ and $L_w$ instead of $\Delta(w\cdot 0)$
and $L(w\cdot 0)$, respectively. Also, for $w\in W$, we denote by
$P_w$ the indecomposable projective cover of $L_w$ in $\mathcal{O}_0$.

The module $\Delta_e$, where $e\in W$ is the identity, is projective
in $\mathcal{O}_0$. The corresponding $L_e$ is the trivial 
$\mathfrak{g}$-module. If we, as usual, denote by $w_0$ the longest
element in $W$, then we have $\Delta_{w_0}=L_{w_0}$ and this is the
only simple Verma module in $\mathcal{O}_0$.

The problem to determine the answer to Kostant's problem for
simple highest weight module $L_w$, where $w\in W$, is 
still wide open in the general case. In what follows, we will 
give a review of known results on this problem which includes 
some progress made in the recent years. For simplicity,
we will say that $L_w$ is {\em Kostant positive} provided
that the solution to Kostant's problem for $L_w$ is positive,
that is, the natural embedding
\begin{displaymath}
U(\mathfrak{g})/\mathrm{Ann}_{U(\mathfrak{g})}(L_w)
\subset
\mathcal{L}(L_w,L_w)
\end{displaymath}
is an isomorphism, and {\em Kostant negative} otherwise. 
Note that, from Subsection~\ref{s2.5} it follows that 
$L_e$ is Kostant positive. Also, from Theorem~\ref{josephthm}
it follows that $L_{w_0}$ is Kostant positive.

We will start or review with a classical results of Gabber and Joseph.

\subsection{Gabber-Joseph's Theorem}\label{s2.8}

Let $W'$ be a parabolic subgroup of $W$ and $w'_0$ the longest
element in $W'$. The following result is due to 
Gabber and Joseph, see \cite[Theorem~4.4]{GJ}

\begin{theorem}[Gabber-Joseph, \cite{GJ}]\label{gabbertheorem}
The module $L_{w'_0w_0}$ is Kostant positive. 
\end{theorem}

To explain this result (which we will do in Subsection~\ref{s7.2}), 
we need to go deeper into the 
properties of $\mathcal{O}$ and explore the connection
between Kostant's problem and projective functors on $\mathcal{O}$.

\section{Projective functors}\label{s3}

\subsection{Definition}\label{s3.1}

For each finite dimensional $\mathfrak{g}$-module $V$, we have the 
endofunctor $V\otimes_\mathbb{C}{}_-$ of the category of all
$\mathfrak{g}$-modules. This endofunctor preserves category 
$\mathcal{O}$. By definition, see \cite{BG}, a {\em projective 
endofunctor} of $\mathcal{O}$ is a functor isomorphic to a direct 
sum of a finite set of indecomposable summands of 
the endofunctors of $\mathcal{O}$ of the form $V\otimes_\mathbb{C}{}_-$.
It is important to point out that projective functors are exact 
and closed under taking both left and right adjoint.

We define the bicategory $\mathscr{P}$ of {\em projective functors}
as follows:
\begin{itemize}
\item objects of $\mathscr{P}$ are indecomposable blocks of $\mathcal{O}$;
\item $1$-morphisms of $\mathscr{P}$ are projective functors between
the corresponding blocks;
\item $2$-morphisms of $\mathscr{P}$ are natural transformations of functors.
\end{itemize}
In particular, we have the monoidal category 
$\mathscr{P}_0:=\mathscr{P}(\mathcal{O}_0,\mathcal{O}_0)$ 
of projective endofunctors of the principal block $\mathcal{O}_0$.

By definition, we have the {\em defining birepresentation} of
the bicategory $\mathscr{P}$ given by its action on blocks
of category $\mathcal{O}$. In particular, the monoidal category 
$\mathscr{P}_0$ acts on the category $\mathcal{O}_0$.

\subsection{Classification}\label{s3.2}

Indecomposable projective functors are classified in \cite[Theorem~3.3]{BG}.
Since we will mostly work with the principal block $\mathcal{O}_0$,
we only recall the classification of indecomposable projective functors
in $\mathscr{P}_0$. The claim of \cite[Theorem~3.3]{BG} is that 
there is a bijection between the elements of $W$
and isomorphism classes of indecomposable projective endofunctors of
$\mathcal{O}_0$. For $w\in W$, there is a unique indecomposable 
projective endofunctor of $\mathcal{O}_0$ such that $\theta_w P_e\cong P_w$.

\subsection{Gradings}\label{s3.3}

Denote by $A$ the opposite of the endomorphism algebra of the direct sum of all
$P_w$, where $w\in W$. The latter direct sum is a projective generator of
$\mathcal{O}_0$. This means, in particular, that $\mathcal{O}_0$ is equivalent 
to the category of finite dimensional $A$-modules. By \cite{So}, the
algebra $A$ is Koszul, in particular, it admits a positive $\mathbb{Z}$-grading
\begin{displaymath}
A=\bigoplus_{i\geq 0}A_0,
\end{displaymath}
where $A_0$ is semi-simple and the positive part 
$\displaystyle \bigoplus_{i> 0}A_0$ coincides with the radical of 
the algebra $A$.

This means that we can consider the category ${}^{\mathbb{Z}}\mathcal{O}_0$
of finite dimensional $\mathbb{Z}$-graded $A$-modules (the morphisms
in this category are homogeneous homomorphisms of degree zero). This is what is 
usually called the {\em graded lift} of $\mathcal{O}_0$. All structural
modules admit graded lifts. For indecomposable modules such lifts are
necessarily unique up to isomorphism and shift of grading. For
$w\in W$, one can fix the graded lift of $L_w$ in degree zero
and then graded lifts of $\Delta_w$ and $P_w$ such that the
top is in degree zero.

For $k\in\mathbb{Z}$, we have the auto-equivalence $\langle k\rangle$
of ${}^{\mathbb{Z}}\mathcal{O}_0$ which shifts the degree by subtracting $k$.

Projective functors also admit graded lifts, so that we have the 
corresponding monoidal category ${}^{\mathbb{Z}}\mathscr{P}_0$
of graded projective functors acting on ${}^{\mathbb{Z}}\mathcal{O}_0$.
We fix the grading lift of $\theta_w$ such that $\theta_w P_e\cong P_w$,
just like in the ungraded case, but now meaning that no shifts of
gradings are involved. We refer to \cite{St} for details.

\subsection{Combinatorics}\label{s3.4}

Let $S$ be the set of simple reflections in $W$
associated with our choice of $\mathbf{R}_+$.
Consider the Hecke algebra $\mathbf{H}$ of 
the Coxeter system $(W,S)$ over
$\mathbb{Z}[v,v^{-1}]$. This is the associative algebra
generated by $H_s$, where $s\in S$, subject to the same
braid relations as the corresponding simple reflections
and, in addition, to the quadratic relations 
$(H_s+v)(H_s-v^{-1})=0$, where $s\in S$. This algebra 
has the {\em standard basis} $\{H_w\,:\, w\in W\}$ and the
{\em Kazhdan-Lusztig basis} $\{\underline{H}_w\,:\,w\in W\}$,
see \cite{So2}.

Fix an isomorphism, $\Phi$, between $\mathbf{H}$ and the Grothendieck group
of ${}^{\mathbb{Z}}\mathcal{O}_0$ such that each $H_w$ corresponds
to $[\Delta_w]$ and such that $\langle 1\rangle$ corresponds to multiplication by $v^{-1}$.
As follows from the Kazhdan-Lusztig Conjecture (which is a theorem,
see \cite{KL,BB,BK,EW}), 
under $\Phi$, the class of $[P_w]$ corresponds to $\underline{H}_w$,
for $w\in W$. Similarly, we have an isomorphism between $\mathbf{H}$
and the split Grothendieck ring of ${}^{\mathbb{Z}}\mathscr{P}_0$
that matches $\underline{H}_w$ with $[\theta_w]$.
In this way, the action of ${}^{\mathbb{Z}}\mathscr{P}_0$ on 
${}^{\mathbb{Z}}\mathcal{O}_0$ gives a right regular $\mathbf{H}$-module.

For $x,y\in W$, we write $\theta_x\geq_R \theta_y$ provided that 
there is $z\in W$ such that $\theta_x$ is isomorphic to a direct summand
of $\theta_z\circ \theta_y$. This is the {\em Kazhdan-Lusztig (KL) 
right pre-order}. Equivalence classes with respect to it are called
{\em right cells}. Similarly one defines the left pre-order and the
left cells using composition on the other side, and 
the two-sided pre-order and the two-sided cells using composition
on both sides. Note the unnaturality of the terminology 
(with respect to left and right) that  stems from the right nature 
of the action of ${}^{\mathbb{Z}}\mathscr{P}_0$ on 
${}^{\mathbb{Z}}\mathcal{O}_0$.

For $x,y\in W$, we set $x\geq_R y$ if and only if $\theta_x\geq_R \theta_y$
and similarly for $L$ and $J$. In this way we can speak of left, right and
two-sided orders and cells on $W$.

\subsection{Annihilators}\label{s3.5}

If we want to study Kostant's problem for simple highest weight modules
$L_w$, where $w\in W$, it is important to understand annihilators of
these modules. The following result is due to Barbasch and Vogan, 
see \cite{BV1,BV2}.

\begin{theorem}[Barbasch-Vogan, \cite{BV1,BV2}]
\label{barbaschvaganthm}
For any $x,y\in W$, we have the containment  
$\mathrm{Ann}_{U(\mathfrak{g})}(L_x)\subset 
\mathrm{Ann}_{U(\mathfrak{g})}(L_y)$ if and only if 
$x\geq_L y$. In particular, we have
$\mathrm{Ann}_{U(\mathfrak{g})}(L_x)= 
\mathrm{Ann}_{U(\mathfrak{g})}(L_y)$
if and only if $x$ and $y$ belong to the same
Kazhdan-Lusztig left cell.
\end{theorem}

\subsection{Twisting functors}\label{s3.6}

Let $s\in W$ be a simple reflection and $\alpha\in\mathbf{R}_-$
the corresponding negative root. Fix a non-zero element 
$Y_\alpha\in\mathfrak{g}_\alpha$ and consider the localization
$U_\alpha$ of $U(\mathfrak{g})$ with respect to 
$\{Y_\alpha^{i}\,:\,i\in\mathbb{Z}_{\geq 0}\}$.
Note that $U_\alpha$ is an algebra and that $U(\mathfrak{g})$ is a 
$U(\mathfrak{g})$-$U(\mathfrak{g})$-subbimodule of the
$U(\mathfrak{g})$-$U(\mathfrak{g})$-bimodule $U_\alpha$.

Let $\phi_s$ denote an inner automorphism of
$\mathfrak{g}$ which swaps each $\mathfrak{g}_\beta$,
where $\beta\in\mathbf{R}$, with $\mathfrak{g}_{s(\beta)}$.
Then we have the $U(\mathfrak{g})$-$U(\mathfrak{g})$-bimodule
$B_s:={}^{\phi_s}(U_\alpha/U(\mathfrak{g}))$ which is obtained
from $U_\alpha/U(\mathfrak{g})$ by twisting the left action of
$U(\mathfrak{g})$ by $\phi_s$.

For a simple reflection $s\in W$, we have the corresponding
{\em twisting functor} $\top_s$ on $\mathcal{O}_0$, which is defined
as tensoring with $B_s$, see \cite{AS}. Note that the twist by
$\phi_s$ is necessary in order to have a functor which preserves 
$\mathcal{O}$. The functors $\{\top_s\,:\,s\in S\}$ satisfy 
braid relations, see \cite{KM}, which allows us to define 
the twisting functors $\top_w$, for each $w\in W$, in the usual
way using reduced expressions. Twisting functors have a number of 
remarkable properties. Here is a brief summary of those which 
are relevant for our goals, see \cite{AS,KM} for details:
\begin{itemize}
\item Twisting functors functorially commute with projective
functors.
\item Twisting functors are right exact and acyclic on Verma modules.
\item The left derived functors of twisting functors are
self-equivalences of $\mathcal{D}^b(\mathcal{O})$.
\item $\top_w \Delta_e\cong \Delta_w$, for each $w\in W$.
\end{itemize}

\subsection{Relevance and potential strategy}\label{s3.7}

Relevance of projective functors for Kostant's problem is two fold.
On the one hand, as we have seen in Theorem~\ref{barbaschvaganthm},
combinatorics of composition of projective functors controls
inclusions between the annihilators of simple highest weight modules
in the principal block $\mathcal{O}_0$.

On the other hand, projective functors also appear in the 
canonical isomorphisms in Subsection~\ref{s1.5}, that is 
on the codomain side of the homomorphism in Kostant's problem.
This suggests a potential strategy for attacking 
Kostant's problem using twisting functors: start from some module
for which the answer is known, apply twisting functors and try to
control both the domain part of Kostant's problem via
annihilators and the codomain part via the
canonical isomorphisms and using the fact that twisting and
projective functors commute.

\section{More recent results on Kostant's problem}\label{s4}

\subsection{Kostant's problem for Verma modules using twisting functors}\label{s4.1}

We can, finally, explain the proof of Theorem~\ref{josephthm}
for all Verma modules in $\mathcal{O}_0$. As we have already mentioned
in Subsection~\ref{s2.4}, the answer to Kostant's problem is positive
for $\Delta_e$. We also know that all $\Delta_w$ have the same
annihilator. It remains to show that the right hand side of 
\eqref{eq2} is the same for all Verma modules. For this we take
$\mathrm{Hom}_\mathfrak{g}(\Delta_e,\Delta_e\otimes_\mathbb{C} V^*)$
and apply $\top_w$ to both arguments of this space. 
Since $\top_w$ is a derived equivalence that
is acyclic on Verma modules and commutes with projective functors,
we get the isomorphic space
$\mathrm{Hom}_\mathfrak{g}(\Delta_w,\Delta_w\otimes_\mathbb{C} V^*)$
using $\top_w\Delta_e=\Delta_w$. That is it.

\subsection{Some further results}\label{s4.2}

Various refinements of the approach via twisting functors were used in
\cite{MS,Ma,KaM,Ka} to obtain a number of results on Kostant's problem for
simple modules in $\mathcal{O}_0$. We start with the following explicit
result by the author, see \cite[Theorem~1]{Ma}:

\begin{theorem}\label{mazthm}
Let $W'$ be a parabolic subgroup of $W$ and $s\in W'$ be a 
simple reflection. Then the module $L_{sw'_0w_0}$ is Kostant positive. 
\end{theorem}

This theorem should be compared to Theorem~\ref{gabbertheorem}.
As promised, we will comment more on Theorem~\ref{gabbertheorem}
later on when we will talk about the parabolic version of
category $\mathcal{O}$, see Subsection~\ref{s7}.

Theorem~\ref{mazthm} was further generalized and extended to a
comparison result by K{\aa}hrtsr{\"o}m, see \cite[Theorem 1.1]{Ka}:

\begin{theorem}[K{\aa}hrtsr{\"o}m, \cite{Ka}]\label{kahrstromthm}
Let $W'$ be a parabolic subgroup of $W$ and $w\in W'$.
Let $\mathfrak{a}$ be the semi-simple Lie algebra whose
Weyl group is $W'$. Then the module $L_w$ is Kostant positive
in the category $\mathcal{O}$ for $\mathfrak{a}$ if and only if
the module $L_{ww'_0w_0}$ is Kostant positive in the category 
$\mathcal{O}$ for $\mathfrak{g}$. 
\end{theorem}

\subsection{Left cell invariance in type $A$}\label{s4.5}

Recall from Subsection~\ref{s3.5} that simple highest weight 
modules in $\mathcal{O}_0$ have the same annihilator in
$U(\mathfrak{g})$ if and only if their indices belong to the
same left cell. As it turned out, in type $A$, the answer 
to Kostant's problem for a simple highest weight module
in $\mathcal{O}_0$ is also constant on a left cell.
The following theorem is due to Catharina Stroppel and the 
author, see \cite[Theorem~V]{MS2}.

\begin{theorem}\label{stroppelmazorchukthm}
Assume that we are in type $A$, that is $W=S_n$. 
Let $x,y\in W$ be two elements in the same 
left cell. Then $L_x$ is Kostant positive 
if and only if $L_y$ is.
\end{theorem}

In type $A$, the cell combinatorics can be described
in terms of the Robinson-Schensted correspondence, 
see \cite[Subsection~3.1]{Sa} for the latter and 
\cite[Section~5]{KL} for the connection between the two.
Recall that the Robinson-Schensted correspondence
associates to any element $w\in S_n$ a pair
$(\alpha(w),\beta(w))$ of standard Young tableaux
of the same shape (and this shape is a partition of $n$). It turns out that
\begin{itemize}
\item $x,y\in S_n$ belong to the same left cell
if and only if $\beta(x)=\beta(y)$;
\item $x,y\in S_n$ belong to the same right cell
if and only if $\alpha(x)=\alpha(y)$;
\item $x,y\in S_n$ belong to the same two-sided cell
if and only if the shape of $\alpha(x)$ coincides
with the shape of $\alpha(y)$.
\end{itemize}

In particular, each left cell contains a unique involution
(involutions in $S_n$ are characterized by the property
that $\alpha(w)=\beta(w)$). Thus, Theorem~\ref{stroppelmazorchukthm}
reduces Kostant's problem for all simple highest weight modules
in $\mathcal{O}_0$ for $\mathfrak{sl}_n$ to Kostant's problem  for 
all simple highest weight modules indexed by involutions.

\subsection{Negative examples}\label{s4.3}

If $\mathfrak{g}$ is of type $B_2$, then $W$ is generated by
two simple reflections, say $s$ and $t$, subject to the relations
$s^2=t^2=e$ and $stst=tsts$. The group $W$ has four left cells,
namely $\{e\}$, $\{w_0\}$, $\{s,ts,sts\}$ and $\{t,st,tst\}$.
By Theorem~\ref{barbaschvaganthm}, the simple modules
$L_s$ and $L_{ts}$ have the same annihilator. At the same time,
it is easy to check, by a direct computation, that 
\begin{displaymath}
\dim\mathrm{Hom}_{\mathfrak{g}}(L_s,\theta_s\theta_t\theta_s L_s)=1
\quad\text{ while }\quad
\dim\mathrm{Hom}_{\mathfrak{g}}(L_{ts},\theta_s\theta_t\theta_s L_{ts})=2.
\end{displaymath}
Taking the canonical isomorphisms \eqref{eq2} into account, 
this implies that at most one of these two modules can be 
Kostant positive. Since $L_s$ is Kostant positive by
Theorem~\ref{mazthm}, it follows that $L_{ts}$ is Kostant negative.
This was already observed in \cite[Corollary~9.5]{Jo}, see
also \cite[Subsection~11.5]{MS2}. 

It is much more tricky to construct an example of a Kostant
negative module in type $A$. Such example, which was a fair surprise, 
was first constructed
in \cite[Theorem~12]{MS}. Consider $\mathfrak{g}$ of type $A_3$
with $W$ generated by three simple reflections $r$, $s$ and $t$,
such that $r$ and $t$ commute. Then the module $L_{rt}$ is
Kostant negative. The proof is similar in spirit 
to the proof in the previous paragraph. However, here the module
$L_{rt}$ is not compared to another simple highest weight module,
but rather to the primitive quotient of the projective Verma module
(and we know that all such quotients are Kostant positive) by
the submodule generated by the annihilator of $L_{rt}$.

Since the answer to Kostant's problem is a left cell invariant in type
$A$ by Theorem~\ref{stroppelmazorchukthm}, it follows that the 
module $L_{srt}$ is Kostant negative as well.

\subsection{Kostant's problem and double centralizers}\label{s4.4}

The paper \cite{KaM} makes a deep analysis of the negative example
in type $A_3$ mentioned in the previous subsection. 
The outcome is a module-theoretic criterion for Kostant positivity of 
simple highest weight modules in $\mathcal{O}_0$, in particular,
in terms of a certain double centralizer property.

Assume that we are in type $A$. Then, as explained above,
it is enough to understand Kostant's problem for modules
$L_d$, where $d$ is an involution. So, let $d$ be an involution.
Let $\theta\in \mathscr{P}_0$
and $x\in W$ be such that $L_x$ is a composition subquotient of
$\theta L_d$. Then $x\leq_R d$, see \cite[Lemma~13(a)]{MM}.

Let $\mathtt{R}$ denote the right cell of $d$ and
$\hat{\mathtt{R}}$ the set of all $y\in W$ such that $y\leq_R d$.
Consider the Serre subcategory $\mathcal{O}_0^{\hat{\mathtt{R}}}$
of $\mathcal{O}_0$ generated by all $L_y$, where $y\in \hat{\mathtt{R}}$.
Then, by the above, the action of $\mathscr{P}_0$ on $\mathcal{O}_0$
restricts to $\mathcal{O}_0^{\hat{\mathtt{R}}}$. Furthermore,
from the canonical isomorphisms \eqref{eq2} it follows that 
all information necessary for determination of Kostant 
positivity or negativity of $L_w$ can be found already inside
this category $\mathcal{O}_0^{\hat{\mathtt{R}}}$.
In other words, we can replace $\mathcal{O}_0$
by $\mathcal{O}_0^{\hat{\mathtt{R}}}$.

The category $\mathcal{O}_0^{\hat{\mathtt{R}}}$ was studied in detail
in \cite{MS2,Ma2}. Similarly to $\mathcal{O}_0$, it turns out that 
$\mathcal{O}_0^{\hat{\mathtt{R}}}$ has non-trivial projective-injective modules.
In fact, the indecomposable projective cover
$P_x^{\hat{\mathtt{R}}}$ of $L_x$ in $\mathcal{O}_0^{\hat{\mathtt{R}}}$, 
where $x\in \hat{\mathtt{R}}$, is injective if and only if 
$x\in\mathtt{R}$.

The projective-injective module $P_x^{\hat{\mathtt{R}}}$, for
$x\in\mathtt{R}$, can be also explicitly constructed. 
For $x\in\mathtt{R}$, we have $P_x^{\hat{\mathtt{R}}}\cong\theta_x L_d$.
In particular, the module $\theta_d L_d$ is an indecomposable 
projective-injective module in $\mathcal{O}_0^{\hat{\mathtt{R}}}$
with simple top and simple socle, both isomorphic to $L_d$.
Since $\theta_d$ is self-adjoint, we have
\begin{displaymath}
\mathrm{Hom}_{\mathfrak{g}}(P_e,\theta_d L_d)\cong
\mathrm{Hom}_{\mathfrak{g}}(\theta_d P_e,L_d)\cong
\mathrm{Hom}_{\mathfrak{g}}(P_d,L_d)\cong \mathbb{C}.
\end{displaymath}
This means that there is a unique, up to scalar, non-zero map from
$P_e$ to $\theta_d L_d$. Let $Q_d$ denote the image of such map.
The main result of \cite{KaM} is the following theorem,
see \cite[Theorem~5]{KaM}.

\begin{theorem}\label{kahrstrommazorchukthm}
The module $L_d$ is Kostant positive if and only if each simple
subquotient in the socle of $(\theta_d L_d)/Q_d$ is of the form
$L_x$, for $x\in\mathtt{R}$.
\end{theorem}

An equivalent reformulation of the socle condition in 
Theorem~\ref{kahrstrommazorchukthm} is as follows:
the module $Q_d$ admits a two-step coresolution
\begin{equation}\label{eq3}
0\to  Q_d\to \theta_d L_d\to N,
\end{equation}
where $N$ is projective-injective (note that $\theta_d L_d$ is projective-injective
as well).

In \cite[Proposition~2]{Ma2}, it was shown that 
$Q_d$ is isomorphic to $P_e^{\hat{\mathtt{R}}}$.
Since all projective modules in $\mathcal{O}_0^{\hat{\mathtt{R}}}$
can be obtained, up to isomorphism, by applying projective
functors to $P_e^{\hat{\mathtt{R}}}$, and projective functors are
exact and preserve both projective and injective objects, it follows
that, under the existence of \eqref{eq3},
any projective object $P\in \mathcal{O}_0^{\hat{\mathtt{R}}}$
has a two-step coresolution
\begin{displaymath}
0\to  P\to N_0\to N_1,
\end{displaymath}
in which both $N_0$ and $N_1$ are projective-injective. This means that the
associative algebra of $\mathcal{O}_0^{\hat{\mathtt{R}}}$ has dominance dimension
at least two with respect to projective injective objects which, in turn, is 
equivalent to the double-centralizer property for this algebra 
on projective-injective objects, see \cite{KSX} for details. Therefore Kostant 
positivity  of $L_d$ is equivalent to this double centralizer property,
see \cite[Corollary~12]{Ma2}.

\subsection{Explicit results in small ranks}\label{s4.6}

Theorem~\ref{kahrstrommazorchukthm} allows one to approach Kostant's
problem for simple highest weight modules in $\mathcal{O}_0$ by 
direct computations (in particular, using computer), especially
in small rank cases.

It is easy to check that in types $A_1$ and $A_2$ all simple 
highest weight modules in $\mathcal{O}_0$ are Kostant positive.
In type $A_3$ we have two only two Kostant negative modules, namely
$L(rt)$ and $L(srt)$, as already mentioned in Subsection~\ref{s4.3}.
In type $A_4$, a complete answer is given in \cite[Proposition~25]{KaM}.
The type $A_5$ was started in \cite{KaM}, further considered
in \cite{Ka} and, finally, completed in \cite[Subsection~10.1]{KMM}.
Additionally, \cite[Subsection~10.2]{KMM} completely answers
Kostant's problem in types $B_3$, $B_4$ and $D_4$.
These computations turned out to be very useful 
for some of the results which will be discussed later.

\section{Kostant's problem and birepresentation theory}\label{s5}

\subsection{K{\aa}hrstr{\"o}m's conjecture}\label{s5.1}

In March 2019, the author received an email from Johan K{\aa}hrstr{\"o}m
that contained the following conjecture, which was based on a comprehensive
analysis of all known results, in particular, all known small rank cases:

\begin{conjecture}[K{\aa}hrstr{\"o}m, 2019] \label{khconjecture}
For an involution $d\in S_n$, the following conditions are equivalent:
\begin{enumerate}[(a)]
\item\label{khconjecture.1} The module $L_d$ is Kostant positive.
\item\label{khconjecture.2} For all $x,y\in S_n$ such that $x\neq y$
and both $\theta_x L_d\neq 0$ and $\theta_y L_d\neq 0$, we have
$\theta_x L_d\not\cong \theta_y L_d$.
\item\label{khconjecture.3} For all $x,y\in S_n$ such that $x\neq y$
and both $\theta_x L_d\neq 0$ and $\theta_y L_d\neq 0$, we have
$[\theta_x L_d]\not=[\theta_y L_d]$ in the Grothendieck group of  
${}^\mathbb{Z}\mathcal{O}_0$.
\item\label{khconjecture.4} For all $x,y\in S_n$ such that $x\neq y$
and both $\theta_x L_d\neq 0$ and $\theta_y L_d\neq 0$, we have
$[\theta_x L_d]\not=[\theta_y L_d]$ in the Grothendieck group of
$\mathcal{O}_0$.
\end{enumerate}
\end{conjecture}

Note that, due to the left cell invariance in type $A$, 
this conjecture effectively covers the whole of $S_n$,
not only the involutions. This conjecture provided a very interesting
new insight into the inner structures that control the behavior
of the answer to Kostant's problem. It motivated various instances of new 
research done in this area and led to substantial progress
and a number of new results. Basically, all the results which will
be described in the remainder of this paper are either directly or
indirectly motivated by this conjecture.

\subsection{Birepresentations of projective functors}\label{s5.2}

As we have already mentioned above, the monoidal category $\mathscr{P}_0$
acts on $\mathcal{O}_0$. Fix some element $w\in W$ and consider the 
corresponding object $L_w\in \mathcal{O}_0$. There are two natural ways
to construct a birepresentation of $\mathscr{P}_0$ starting from
$L_w$.

{\em First construction.} Let $\mathscr{I}_w$ denote the annihilator
of $L_w$ in $\mathscr{P}_0$. This is a left monoidal ideal in 
$\mathscr{P}_0$. The corresponding quotient $\mathscr{P}_0/\mathscr{I}_w$
has the natural (left) action of $\mathscr{P}_0$ induced from the regular
(left) action of $\mathscr{P}_0$ on itself. We denote this
birepresentation of $\mathscr{P}_0$ by $\mathbf{M}$.

{\em Second constriction.} There is the obvious action of $\mathscr{P}_0$
on the additive closure (inside $\mathcal{O}_0$) of all objects of
the form  $\theta L_w$, where $\theta\in \mathscr{P}_0$. We denote this
birepresentation of $\mathscr{P}_0$ by $\mathbf{N}$.

There is the obvious morphism of birepresentation from $\mathbf{M}$
to $\mathbf{N}$ which sends $\theta\in \mathbf{M}$ to
$\theta L_w\in \mathbf{N}$. The following observation follows
from \cite[Subsection~8.3]{KMM} by comparing 
\cite[Proposition~8.11]{KMM} with  \cite[Formula~(39)]{KMM}:

\begin{proposition}\label{kmmprop}
For $w\in W$, the module $L_w$ is Kostant positive if and only if
the natural morphism of birepresentation of $\mathscr{P}_0$ 
from $\mathbf{M}$ to $\mathbf{N}$ is an equivalence.
\end{proposition}

An interesting consequence of this observation is that one can 
now formulate and study Kostant's problem  in the general setup
of birepresentations of finitary bicategories, see \cite[Subsection~7.2]{MMM}.

\subsection{Indecomposability conjecture}\label{s5.3}

K{\aa}hrstr{\"o}m's conjecture suggests that it is necessary to understand
the modules of the form $\theta_x L_y$, where $x,y\in W$. Such modules
played important roles in various other investigations, see \cite{Ma2,KiM}.
In particular,  the paper \cite{KiM} formulates the following 
indecomposability conjecture,  see \cite[Conjecture~2]{KiM}:

\begin{conjecture}\label{kimconjecture}
For $x,y\in S_n$, the module $\theta_x L_y$ is either zero or indecomposable. 
\end{conjecture}

\subsection{K{\aa}hrstr{\"o}m's conjecture vs indecomposability conjecture}\label{s5.4}

One of the main applications of Proposition~\ref{kmmprop} in 
\cite{KMM} is the following result, see \cite[Theorem~8.16]{KMM}:

\begin{theorem}\label{kmmthm}
Let $d\in S_n$ be an involution. Then Conjecture~\ref{khconjecture}\eqref{khconjecture.1}
is equivalent to the conjunction of Conjecture~\ref{khconjecture}\eqref{khconjecture.2}
with the validity of Conjecture~\ref{kimconjecture} for $y=d$ and for all $x\in S_n$.
\end{theorem}

In fact, \cite[Theorem~8.16]{KMM} is formulated for arbitrary type, but under
the assumption that $d$ is a so-called {\em distinguished} (or {\em Duflo}) involution,
and not just an involution. Distinguished involutions are special elements in
left (right) cells. In fact, each left (right) cell contains a unique 
distinguished involution, see \cite[Proposition~17]{MM}.

\subsection{Equivalences of categories using Kostant's problem}\label{s5.5}

Proposition~\ref{kmmprop} has the following consequence: given two 
$x,y\in W$ from the same left cell and such that both 
$L_x$ and $L_y$ are Kostant positive, it follows that the additive
closures of $\mathscr{P}_0 L_x$ and $\mathscr{P}_0 L_y$ are equivalent.
Indeed, both of them are equivalent to the quotient of 
$\mathscr{P}_0$ by the (common) annihilator of these two modules
in $\mathscr{P}_0$.

There is a small subtlety here as the fact that $x$ and $y$ belong to
the same left cell only implies that the annihilators of $L_x$ and $L_y$
in $U(\mathfrak{g})$ coincide. However, it is not very difficult to
deduce from this that the annihilators of $L_x$ and $L_y$
in $\mathscr{P}_0$ coincide as well, see \cite[Proposition~8.11]{KMM}.

The idea to compare two abstract categories using the category of 
projective functors which acts on both of these abstract categories 
goes back to \cite{MiSo},
see \cite[Theorem~3.1]{MiSo}. In turn, \cite{MiSo} attributes this 
idea to the paper \cite{BG} where projective functors were defined 
and classified. The classification in \cite{BG} is achieved by 
describing the category of projective functors in terms of a
certain category of Harish-Chandra bimodules. In any case, the input
of \cite[Theorem~3.1]{MiSo} is a certain Kostant positive 
$\mathfrak{g}$-module $M$ with a fixed annihilator. The output is
an equivalence similar to that given by Proposition~\ref{kmmprop}.
Such equivalences turned out to be very useful, for example, 
to describe  the structure of induced modules over Lie algebras, 
see \cite{KM0,MS2} and references therein. Such an equivalence
allows one to connect the category one wants to understand with some
well-known category, for example, with category $\mathcal{O}$.
After that one can track the modules one wants to understand through
the equivalence and read off the structural information from that
in, say, category $\mathcal{O}$. In particular, this shows that 
knowing the answer to Kostant's problem for particular modules could 
be very useful.

\section{Kostant's problem for fully commutative permutations}\label{s6}

\subsection{Fully commutative permutations}\label{s6.1}

Consider $S_n$ as a Coxeter group in the usual way, with simple reflections
given by the elementary transpositions $s_i=(i,i+1)$, where $i=1,2,\dots,n-1$. Given
$w\in S_n$, we might have several {\em reduced expressions} for $w$, that
is shortest possible factorizations of $w$ into a product of simple
reflections. A classical result of Matsumoto, see \cite{Mat}, asserts that
any two reduced expressions of $w$ can be transformed one into other
in a finite number of steps each of which uses only some braid relations,
that is the relations of the form
\begin{displaymath}
s_is_j=s_js_i,\quad |i-j|>1;\qquad
s_is_{i+1}s_i=s_{i+1}s_is_{i+1}.
\end{displaymath}
An element $w$ is called {\em fully commutative} provided that 
any two reduced expressions of $w$ can be transformed one into other
in a finite number of steps each of which uses only the commutativity 
relations $s_is_j=s_js_i$, where $|i-j|>1$.

An alternative characterization of fully commutative permutations can
be given in terms of the Robinson-Schensted correspondence
(the latter correspondence was already mentioned in Subsection~\ref{s4.5}).
An element $w\in S_n$ is fully commutative if and only if the shape of
its Robinson-Schensted correspondent has at most two rows.

\subsection{Fully commutative part of $\mathcal{O}_0$}\label{s6.2}

The Serre subcategory of $\mathcal{O}_0$ for $\mathfrak{sl}_n$ 
generated by all $L_w$ such that $w\in S_n$ is fully commutative
admits a very nice and useful diagrammatic description developed
in \cite{BS1,BS2}. One particular  consequences of this description
is the following statement, see \cite[Theorem~4.11(iii)]{BS1}:

\begin{theorem}[Brundan-Stroppel, \cite{BS1}]\label{bstheorem}
Let $y\in S_n$  be a fully commutative element and $x\in S_n$ be such that 
$\theta_x L_y\neq 0$. Then $x$  is fully commutative as well and
$\theta_x L_y$ has simple  top, in particular, $\theta_x L_y$ is indecomposable.
\end{theorem}

Combining this with Theorem~\ref{kmmthm}, we obtain that,  for
a fully commutative involution $d\in S_n$, the first two assertions
in Conjecture~\ref{khconjecture} are equivalent. This was
a staring observation for the project that resulted in the paper
\cite{MMM}. In order to describe the results of this paper,
we need some definitions.

\subsection{Special fully commutative involutions}\label{s6.3}

Fix two elements $i<j$ from the set $\{1,2,\dots,n\}$ such that $i+j+1\leq n$.
Denote by $\sigma_{i,j}$ the following involution in $S_n$:

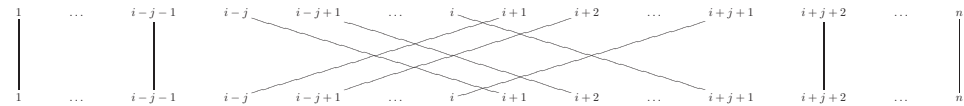
\resizebox{\textwidth}{!}{
$
\xymatrix@R=20mm@C=12mm{
1\ar@{-}[d]&\dots&i-j-1\ar@{-}[d]&i-j\ar@{-}[drrrr]&i-j+1\ar@{-}[drrrr]&\dots
&i\ar@{-}[drrrr]&i+1\ar@{-}[dllll]&i+2\ar@{-}[dllll]&\dots&i+j+1\ar@{-}[dllll]&
i+j+2\ar@{-}[d]&\dots&n\ar@{-}[d]\\
1&\dots&i-j-1&i-j&i-j+1&\dots&i&i+1&i+2&\dots&i+j+1&i+j+2&\dots&n
}
$
}

This element is, in fact, fully commutative and will be called {\em special}.
The subset $\{i-j,i-j+1,\dots,i+j+1\}$ of $\{1,2,\dots,n\}$ on which 
$\sigma_{i,j}$ operates non-trivially is called the {\em support} of 
$\sigma_{i,j}$. Two special elements are called {\em distant} provided
that there is 
an integer $k$ which is strictly greater than all elements in one of
these supports and which is strictly smaller than all elements in the 
other support. In other words, the two supports are separated by
some integer inside $\mathbb{Z}$, considered as a poset in the usual way.

\subsection{Answer to Kostant's problem for fully commutative permutations}\label{s6.4}

The main result  of \cite{MMM} is the following statement,
see \cite[Theorem~5]{MMM}:

\begin{theorem}\label{mmmthm}
K{\aa}hrstr{\"o}m's Conjecture~\ref{khconjecture} is true for all fully commutative 
involutions in $S_n$. Moreover, a fully commutative involution $d\in S_n$
is Kostant positive if and only if $d$ is a product of pairwise distant
special involutions.
\end{theorem}

The strategy of the proof of this theorem is as follows: 
as already mentioned above, at this stage we know that, for a fully 
commutative involution $d\in S_n$, the first two assertions
in Conjecture~\ref{khconjecture} are equivalent. We also have the 
obvious implications
\begin{displaymath}
\text{Conjecture~\ref{khconjecture}\eqref{khconjecture.4}}\Rightarrow
\text{Conjecture~\ref{khconjecture}\eqref{khconjecture.3}}\Rightarrow
\text{Conjecture~\ref{khconjecture}\eqref{khconjecture.2}}.
\end{displaymath}
Now the proof splits into two parts. In the first part, for potentially
Kostant positive elements, it is shown that any fully commutative involution 
$d\in S_n$ which is a product of pairwise distant
special involutions has the property described by 
Conjecture~\ref{khconjecture}\eqref{khconjecture.4}.
In the second  part, for potentially
Kostant negative elements, it is shown that any fully commutative involution 
$d\in S_n$ which is not a product of pairwise distant
special involutions has the property described by the negation of  
Conjecture~\ref{khconjecture}\eqref{khconjecture.2}.

Both parts have a number of case-by-case verifications and crucially use
the diagrammatic description from \cite{BS1,BS2} and related combinatorics
of the Temperley-Lieb algebra. The main insight provided by Conjecture~\ref{khconjecture}
for the proof
is that it allows to reformulate a ring-theoretic property involved in
Kostant's problem in terms of solutions of some combinatorial equations for some
diagrammatic semigroups. And, in this particular case, everything could be solved
by hands without any involvement of computers.

\subsection{Combinatorial consequences}\label{s6.5}

Theorem~\ref{mmmthm} is explicit enough to allow for enumeration of
Kostant positive fully commutative elements. Here is the summary
of the corresponding results, see \cite[Section~6]{MMM} for details:
\begin{itemize}
\item The number of fully commutative involutions in $S_n$
equals $\binom{n}{\lfloor\frac{n}{2}\rfloor}$.
\item The number of Kostant positive fully commutative involutions in $S_n$
equals the $n$-th Fibonacci number.
\item The number of fully commutative elements in $S_n$
equals the $n$-th Catalan number $\frac{1}{n+1}\binom{2n}{n}$.
\item The number of Kostant positive fully commutative elements in $S_n$
equals 
\begin{displaymath}
\sum_{a=0}^{\lfloor\frac{n}{2}\rfloor} \binom{n-a}{a}
\frac{n!(n-2a+1)!}{a!(n-2a)!(n-a+1)!}.
\end{displaymath} 
\end{itemize}
As a consequence, if we let $n$ go to infinity, then, asymptotically,
both almost all fully commutative involutions and almost all
fully commutative elements are Kostant negative.

There is an interesting twist in this tale, though. 
There is a natural refinement of the above statistics. 
Recall that the Robinson-Schensted correspondent of a fully commutative
element has the shape of a partition with at most two parts. The size of the 
second (smaller) part,
which is an integer between $0$ and $\lfloor\frac{n}{2}\rfloor$, coincides
with the value of  Lusztig's $\mathbf{a}$-function from \cite{Lu}
on our element,  see \cite[Corollary~6.6]{MT} where this connection
is explained in detail.

We can now consider fully commutative involutions or fully commutative elements with a fixed 
$\mathbf{a}$-value $a$. Then we have the following enumeration,
see \cite[Section~6]{MMM}, for details:
\begin{itemize}
\item The number of fully commutative involutions in $S_n$
with a fixed $\mathbf{a}$-value $a$
equals $\frac{n!(n-2a+1)!}{a!(n-2a)!(n-a+1)!}$.
\item The number of Kostant positive fully commutative involutions in $S_n$
with a fixed $\mathbf{a}$-value $a$
equals $\binom{n-a}{a}$.
\item The number of fully commutative elements in $S_n$
with a fixed $\mathbf{a}$-value $a$
equals $\left(\frac{n!(n-2a+1)!}{a!(n-2a)!(n-a+1)!}\right)^2$.
\item The number of Kostant positive fully commutative elements in $S_n$
with a fixed $\mathbf{a}$-value $a$
equals $\binom{n-a}{a}
\frac{n!(n-2a+1)!}{a!(n-2a)!(n-a+1)!}$.
\end{itemize}
Note that the factor $\frac{n!(n-2a+1)!}{a!(n-2a)!(n-a+1)!}$ is just the
number of standard Young tableaux of shape $(n-a,a)$ and is computed by
the Hook Formula,  see \cite[Theorem~3.10.2]{Sa}.

As a consequence, if we fix $a$ and let $n$ go to infinity, then, asymptotically,
both almost all fully commutative involutions with our fixed $\mathbf{a}$-value $a$ and almost all
fully commutative elements with our fixed $\mathbf{a}$-value $a$ are Kostant positive.

\section{Kostant's problem for parabolic Verma modules}\label{s7}

\subsection{Parabolic category $\mathcal{O}$}\label{s7.1}

Let $\mathfrak{p}$ be a Lie subalgebra of $\mathfrak{g}$ that 
contains the Borel subalgebra $\mathfrak{h}\oplus \mathfrak{n}_+$.
The algebra $\mathfrak{p}$ is called a {\em parabolic subalgebra}
and it is uniquely determined by a subset 
$\pi_\mathfrak{p}$ of simple roots. Namely, $\mathfrak{p}$
is generated by the Borel subalgebra together with all
$\mathfrak{g}_{-\alpha}$, where $\alpha\in \pi_\mathfrak{p}$.

The parabolic subcategory $\mathcal{O}^\mathfrak{p}$ of 
$\mathcal{O}$ is defined as the full subcategory of $\mathcal{O}$
consisting of all modules the action of $\mathfrak{p}$ on which 
is locally finite, see \cite{RC}. The category $\mathcal{O}^\mathfrak{p}$
is a Serre subcategory of $\mathcal{O}$ and it inherits a block
decomposition from that of $\mathcal{O}$. Let $W(\mathfrak{p})$
be the parabolic subgroup of $W$ corresponding to $\pi_\mathfrak{p}$.
Denote by $W_\mathfrak{p}^\text{short}$ the set of all shortest 
coset representatives in the cosets from
$W(\mathfrak{p})\setminus W$. Then $\mathcal{O}^\mathfrak{p}_0$
is the Serre subcategory of $\mathcal{O}_0$ generated by all
$L_w$, where $w\in W_\mathfrak{p}^\text{short}$.

The natural inclusion of $\mathcal{O}^\mathfrak{p}$ into $\mathcal{O}$
is exact and hence has a left adjoint, denoted $\mathrm{Z}^\mathfrak{p}$
and called the {\em Zuckerman functor}. For $M\in\mathcal{O}$, the module
$\mathrm{Z}^\mathfrak{p} M$ is the maximal quotient of $M$ that 
belongs to $\mathcal{O}^\mathfrak{p}$. Being left adjoint to an exact functor,
$\mathrm{Z}^\mathfrak{p}$ sends projective objects to projective objects.
In particular, the indecomposable projectives in $\mathcal{O}^\mathfrak{p}_0$
are $P_w^\mathfrak{p}:=\mathrm{Z}^\mathfrak{p}P_w$, 
where $w\in W_\mathfrak{p}^\text{short}$.

\subsection{Parabolic Verma modules}\label{s7.2}

The modules $\Delta_w^\mathfrak{p}:=\mathrm{Z}^\mathfrak{p}\Delta_w$, 
where $w\in W_\mathfrak{p}^\text{short}$, are called {\em parabolic Verma modules}.
Recall that, by Theorem~\ref{josephthm}, all Verma modules are
Kostant positive. At the same time, we know that the answer to 
Kostant's problem is not inherited by quotients, in general. 
Therefore the problem to determine the answer to Kostant's problem for 
parabolic Verma modules is not trivial.

There are two cases for which the answer can be deduced from the 
results mentioned above. In the first case, we note that 
the module $\Delta_e^\mathfrak{p}$ is a
quotient of the projective Verma module $\Delta_e$ and hence it
is Kostant positive by Subsection~\ref{s2.5}.

Let $w_0^\mathfrak{p}$ denote the longest element in 
$W(\mathfrak{p})$. Then the element $w_0^\mathfrak{p}w_0$ is the
minimum element in $W_\mathfrak{p}^\text{short}$ with respect to the
Bruhat order. In particular, the module $\Delta_{w_0^\mathfrak{p}w_0}^\mathfrak{p}$
is simple and hence coincides with $L_{w_0^\mathfrak{p}w_0}$.
The fact that the latter module is Kostant positive is exactly the
claim of Theorem~\ref{gabbertheorem}. 

At this stage we can even briefly explain why Theorem~\ref{gabbertheorem} is true.
The point is that, up to homological shift, the module 
$L_{w_0^\mathfrak{p}w_0}$ is the image of $\Delta_e^{\mathfrak{q}}$,
for some other parabolic subalgebra $\mathfrak{q}$, under the derived
twisting functor $\top_{w_0}$.  Since the latter functor preserves the
annihilator and commutes with projective functors, it follows
that it preserves both sides of the equation involved in  Kostant's problem.
Therefore Kostant positivity of $L_{w_0^\mathfrak{p}w_0}$ follows
from Kostant positivity of $\Delta_e^{\mathfrak{q}}$, which holds by the previous paragraph.
The fact that the derived $\top_{w_0}$ sends $\Delta_e^{\mathfrak{q}}$
to $L_{w_0^\mathfrak{p}w_0}$, up to homological shift, 
has another major consequence, namely, that
the category $\mathcal{O}_0^\mathfrak{p}$ is {\em Ringel dual} to 
$\mathcal{O}_0^\mathfrak{q}$, see \cite{So3}.

Kostant's problem for parabolic Verma modules is open, in the general case. 
Below we describe the solutions to this problem in two special cases:
for a minimal parabolic and for a maximal parabolic in type $A$.
In many cases, these solutions show that, up to the action of twisting 
functors, the above  two Kostant positive parabolic Verma modules are 
``the only ones that are Kostant positive''.

\subsection{The case of a minimal parabolic: preliminaries}\label{s7.3}

Fir $k\in\{1,2,\dots,n-1\}$ and consider the parabolic subalgebra 
$\mathfrak{p}_k$ in $\mathfrak{sl}_n$ such that $W(\mathfrak{p}_k)=\{e,(k,k+1)\}$.
Note that $W\cong S_n$.

For $1\leq i< j\leq n$, we denote by $\tau^{n,k}_{i,j}$ the 
element in $S_n$ which is uniquely determined by the following conditions:
\begin{itemize}
\item $\tau^{n,k}_{i,j}(i)=k$;
\item $\tau^{n,k}_{i,j}(j)=k+1$;
\item for all $s<t$ in $\underline{n}\setminus\{i,j\}$, we have
$\tau^{n,k}_{i,j}(s)<\tau^{n,k}_{i,j}(t)$.
\end{itemize}
We denote by $\mathcal{X}_{n,k}$ the set consisting of all these 
$\tau^{n,k}_{i,j}$. Here is the example of the set
$\mathcal{X}_{3,1}$:
\begin{displaymath}
\xymatrix@R=3mm@C=3mm{1\ar@{-}[d]&2\ar@{-}[d]&3\ar@{-}[d]\\1&2&3},\qquad 
\xymatrix@R=3mm@C=3mm{1\ar@{-}[d]&2\ar@{-}[dr]&3\ar@{-}[dl]\\1&2&3},\qquad 
\xymatrix@R=3mm@C=3mm{1\ar@{-}[drr]&2\ar@{-}[dl]&3\ar@{-}[dl]\\1&2&3}.
\end{displaymath}
Further, we split $\mathcal{X}_{n,k}$ into 
a disjoint union of two subsets:
$\mathcal{X}_{n,k}^+$ and $\mathcal{X}_{n,k}^-$,
where
\begin{displaymath}
\mathcal{X}_{n,k}^+=\{\tau^{n,k}_{i,i+1}\,:\,i=1,2,\dots,n-1\}
\quad\text{ and }\quad
\mathcal{X}_{n,k}^-=\mathcal{X}_{n,k}\setminus \mathcal{X}_{n,k}^+.
\end{displaymath}
In the example above, the set $\mathcal{X}_{3,1}^+$ contains the 
leftmost and the rightmost elements while the set 
$\mathcal{X}_{3,1}^-$ contains the middle element.

We denote by $G_k$ the centralizer of $(k,k+1)$ in $S_n$
and by $\hat{G}_k$ the subgroup of $G_k$ that consists of 
all elements in $G_k$ which fix both $k$ and $k+1$.
Note that composition of functions gives rise to a bijection
\begin{displaymath}
\hat{G}_k\times \mathcal{X}_{n,k} \to  (S_n)^{\mathfrak{p}_k}_{\mathrm{short}}.
\end{displaymath}

\subsection{Kostant's problem for parabolic Verma modules
in the case of a minimal parabolic}\label{s7.4}
The following statement  is \cite[Theorem~3]{MSr}.

\begin{theorem}\label{shraddhathm1}
For $w\in W_{\mathfrak{p}_k}^{\mathrm{short}}$, the following assertions are equivalent:
\begin{enumerate}[$($a$)$]
\item\label{shraddhathm1.1} The module 
$\Delta^{\mathfrak{p}_k}_w$ is Kostant positive.
\item\label{shraddhathm1.2} For $x,y\in W$ such that $x\neq y$, 
$\theta_x \Delta^{\mathfrak{p}_k}_w\neq 0$ and
$\theta_y \Delta^{\mathfrak{p}_k}_w\neq 0$, we have $\theta_x \Delta^{\mathfrak{p}_k}_w\not\cong \theta_y \Delta^{\mathfrak{p}_k}_w$ (as ungraded modules).
\item\label{shraddhathm1.3} For all $x,y\in W$ such that $x\neq y$, 
$\theta_x \Delta^{\mathfrak{p}_k}_w\neq 0$ and
$\theta_y \Delta^{\mathfrak{p}_k}_w\neq 0$, we have 
$[\theta_x \Delta^{\mathfrak{p}_k}_w]\neq [\theta_y \Delta^{\mathfrak{p}_k}_w]\langle i\rangle$,
for $i\in\mathbb{Z}$,
in the Grothendieck group of ${}^\mathbb{Z}\mathcal{O}_0$.
\item\label{shraddhathm1.4} For all $x,y\in W$ such that $x\neq y$, 
$\theta_x \Delta^{\mathfrak{p}_k}_w\neq 0$ and
$\theta_y \Delta^{\mathfrak{p}_k}_w\neq 0$, we have 
$[\theta_x \Delta^{\mathfrak{p}_k}_w]\neq [\theta_y \Delta^{\mathfrak{p}_k}_w]$
in the Grothendieck group of $\mathcal{O}_0$.
\item\label{shraddhathm1.5} $w\in \hat{G}_k\circ\mathcal{X}_{n,k}^+$.
\item\label{shraddhathm1.6} The 
annihilator of $\Delta^{\mathfrak{p}_k}_w$ in $U(\mathfrak{g})$
is a primitive ideal.
\end{enumerate}
\end{theorem}

The proof of this theorem in \cite{MSr} is rather non-trivial. The main
reason why it was possible to prove this result is the fact that,
in the case of a minimal parabolic, every parabolic Verma module 
is the cokernel of an inclusion between two usual Verma modules.
The socles of such cokernels (in type $A$) were explicitly described
in \cite{KMM2}. It is the knowledge of these socles which is the crucial
technical tool behind the proof.

\subsection{The case of a maximal parabolic}\label{s7.5}

Now let us consider the case of a maximal parabolic. 
Fix $k\in\{1,2,\dots,n-1\}$ and consider the parabolic subalgebra 
$\mathfrak{q}_k$ in $\mathfrak{sl}_n$ such that $W(\mathfrak{q}_k)$
is the standard embedding of $S_k\times S_{n-k}$ into the Weyl group $S_n$.
In this case, $\mathfrak{q}_k$ is isomorphic to $\mathfrak{sl}_k
\times \mathfrak{sl}_{n-k}$, where $\mathfrak{sl}_k$ is embedded into
the top left corner of $\mathfrak{sl}_n$ and $\mathfrak{sl}_{n-k}$ is embedded into
the bottom right corner of $\mathfrak{sl}_n$.

For $k<\frac{n}{2}$, define the set $\mathcal{Y}_{n,k}$
of $(S_n)_{\mathfrak{q}_k}^{\mathrm{short}}$ as follows:
an element $w\in (S_n)_{\mathfrak{q}_k}^{\mathrm{short}}$
belongs to $\mathcal{Y}_{n,k}$ if and only if there is
$i\in\{1,2,\dots,n\}$ such that $w(i-1+j)=k+j$, for all
$j=1,2,\dots,n-k$. Here is the example of the set
$\mathcal{Y}_{5,2}$:
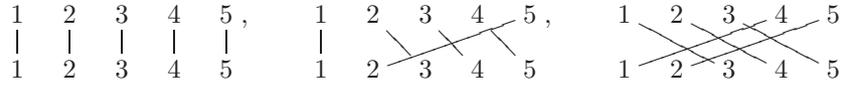
\begin{displaymath}
\xymatrix@R=3mm@C=3mm{1\ar@{-}[d]&2\ar@{-}[d]&3\ar@{-}[d]
&4\ar@{-}[d]&5\ar@{-}[d]\\1&2&3&4&5},\qquad 
\xymatrix@R=3mm@C=3mm{1\ar@{-}[d]&2\ar@{-}[dr]&3\ar@{-}[dr]
&4\ar@{-}[dr]&5\ar@{-}[dlll]\\1&2&3&4&5},\qquad 
\xymatrix@R=3mm@C=3mm{1\ar@{-}[drr]&2\ar@{-}[drr]&3\ar@{-}[drr]
&4\ar@{-}[dlll]&5\ar@{-}[dlll]\\1&2&3&4&5}
\end{displaymath}

For $k>\frac{n}{2}$, define the set $\mathcal{Y}_{n,k}$
of $(S_n)_{\mathfrak{q}_k}^{\mathrm{short}}$ as follows:
an element $w\in (S_n)_{\mathfrak{q}_k}^{\mathrm{short}}$
belongs to $\mathcal{Y}_{n,k}$ if and only if there is
$i\in\{1,2,\dots,n\}$ such that $w(i-1+j)=j$, for all
$j=1,2,\dots,k$. Here is the example of the set
$\mathcal{Y}_{5,3}$:
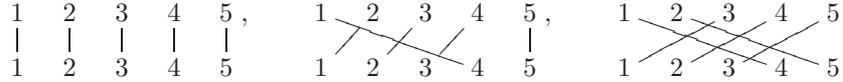
\begin{displaymath}
\xymatrix@R=3mm@C=3mm{1\ar@{-}[d]&2\ar@{-}[d]&3\ar@{-}[d]
&4\ar@{-}[d]&5\ar@{-}[d]\\1&2&3&4&5},\qquad 
\xymatrix@R=3mm@C=3mm{1\ar@{-}[drrr]&2\ar@{-}[dl]&3\ar@{-}[dl]
&4\ar@{-}[dl]&5\ar@{-}[d]\\1&2&3&4&5},\qquad 
\xymatrix@R=3mm@C=3mm{1\ar@{-}[drrr]&2\ar@{-}[drrr]&3\ar@{-}[dll]
&4\ar@{-}[dll]&5\ar@{-}[dll]\\1&2&3&4&5}
\end{displaymath}

\subsection{Kostant's problem for parabolic Verma modules
in the case of a maximal parabolic}\label{s7.6}

The following statement is \cite[Theorem~5]{MSr}:

\begin{theorem}\label{shraddha2}
In the setup above, we have:

\begin{enumerate}[$($a$)$]
\item\label{shraddha2.1} If $k\in\{1,n-1,\frac{n}{2}\}$, then the only 
$w\in W_{\mathfrak{q}_k}^\mathrm{short}$ for which $\Delta^{\mathfrak{q}_k}_w$
is Kostant positive are $w=e$ and $w=w_0^{\mathfrak{q}_k}w_0$.
\item\label{shraddha2.2} If $k\neq\frac{n}{2}$, 
then, for $w\in W_{\mathfrak{q}_k}^\mathrm{short}$, we have that
$\Delta^{\mathfrak{q}_k}_w$
is Kostant positive if and only if $w\in\mathcal{Y}_{n,k}$.
\end{enumerate}
\end{theorem}

In the case of a maximal parabolic, all simple constituents of parabolic
Verma modules are indexed by fully commutative permutations. 
Consequently, the parabolic category $\mathcal{O}_0^{\mathfrak{q}_k}$
admits a diagrammatic description as in \cite{BS1,BS2}. In particular,
both parabolic Verma modules and projective functors admit explicit 
diagrammatic description. This, and Theorem~\ref{mmmthm}, 
are the main ingredients of the proof of \cite[Theorem~5]{MSr}.

\vspace{5mm}

\subsection*{Acknowledgements}

The author is partially supported by the Swedish Research Council.

This paper summarizes the talks given at the
Winter STARS Workshop, Rehovot, Israel, in February 2023;
at the Representation Theory XVIII Conference, Dubrovnik,
Croatia, in June 2023; and at the 14-th Ukraine Algebra Conference 
(via Zoom), Sumy, Ukraine, in July 2023. The author thanks to 
the organizers of these conferences for the invitation
and the opportunity to present these results.

\vspace{2mm}

\noindent
Department of Mathematics, Uppsala University, Box. 480,
SE-75106, Uppsala,\\ SWEDEN, email: {\tt mazor\symbol{64}math.uu.se}

\end{document}